\documentclass[12pt,a4paper]{article}
\usepackage[latin1]{inputenc}
\usepackage[english]{babel}
\usepackage{indentfirst}
\usepackage{amstext,geometry}
\usepackage{amsfonts}
\usepackage{textcomp}
\usepackage{amssymb}
\usepackage{amscd}
\usepackage[dvips]{graphicx}
\usepackage{setspace}
\usepackage{fancyhdr}

\newcommand {\demo}{\hskip -0.6cm{\bf Proof.  }\newline}
\newcommand {\fim}{\nl\text{ }\hfill{$\square$}\vskip 1pc}

\newcommand {\nl}{\newline}
\newcommand {\cl}{\centerline}

\newcommand {\N}{\mathbb{N}}
\newcommand {\Z}{\mathbb{Z}}

\newcommand {\al}{\alpha}

\newcommand {\lb}{\linebreak}

\newcommand {\ovl}{\overline}
\newcommand {\m}[1]{{\widetilde{#1}\,}}
\newcommand{\af}[2]{\vspace{0.3cm}\hskip -0.6cm {\it Claim#1: #2}\vspace{0.3cm}}

\newcommand{\fimaf}{\vspace{0.3cm}}

\newcommand {\prcr}[2]{\mathcal{O}(#1,\al,L_{#2})}
\newcommand {\tp}[2]{\mathcal{T}(#1,\al,L_{#2})}
\newcommand{\funcao}[5]{\begin{array}{rrcl}
#1:&\!\!\!#2 & \rightarrow & #3 \\
  &\!\!\! #4 & \mapsto & #5
\end{array}}
\newcommand{\dfuncao}[7]{\begin{array}{lrcl}
#1:&\!\!\!#2 & \rightarrow & #3 \\
  &\!\!\! #4 & \mapsto & #5\\
  &\!\!\! #6 & \mapsto & #7
\end{array}}

\newtheorem{teorema}{Theorem}[section]
\newtheorem{lema}[teorema]{Lemma}
\newtheorem{corolario}[teorema]{Corollary}
\newtheorem{definicao}[teorema]{Definition}
\newtheorem{proposicao}[teorema]{Proposition}
\begin{document}
\doublespace

\cl{\large\bf The Crossed Product by a Partial Endomorphism}
\cl{\large\bf  and the Covariance Algebra}
\vspace{1cm}
\cl{D. Royer\footnote{Supported by Cnpq}}
\vspace{1cm}

{\abstract Given a local homeomorphism $\sigma:U\rightarrow X$ where $U\subseteq X$ is clopen and $X$ is a
compact and Hausdorff topological space, we obtain the possible transfer operators $L_\rho$ which may occur for
$\al:C(X)\rightarrow C(U)$ given by $\al(f)=f\circ\sigma$. We obtain examples of partial dynamical systems $(X_A,\sigma_A)$
such that the construction of the covariance algebra $C^*(X_A,\sigma_A)$ and the crossed product by partial endomorphism
$\prcr{X_A}{}$ associated to this system are not equivalent, in the sense that there does not exists invertible
function $\rho\in C(U)$ such that $\prcr{X_A}{\rho}=C^*(X_A,\sigma)$.

\section{Introduction}

We start with a summary of the construction of the crossed product by a partial endomorphism. Details
may be seen in \cite{prodcruz}.
A partial $C^*$-dynamical system $(A,\al,L)$ consists of a (closed) ideal $I$ of a $C^*$-algebra $A$, a
idempotent self-adjoint ideal $J$ of $I$ (not necessarily closed), a *-homomorphism $\al:A\rightarrow M(I)$
and a linear posisitve map (which preserves *) $L:J\rightarrow A$ such that $L(a\al(b))=L(a)b$ for each
$a\in J$ and $b\in A$. The map $L$ is called {\it transfer operator}. Define in $J$ an inner product (which may be degenerated) by $(x,y)=L(x^*y)$. Then we obtain a
inner product $\langle,\rangle$ in the quotient $J_0=J / \{x\in J:L(x^*x)=0\}$ defined by
$\langle\m{x},\m{y}\rangle=L(x^*y)$, which induces a norm $\|\|$. Define $M=\ovl{J_0}^{\|\|}$, which is a
right Hilbert $A$-module and also a left $A$-module, where the left multiplication is defined by the
*-homomorphism $\varphi:A\rightarrow L(M)$ (the adjointable operators in $M$), where $\varphi(a)(\m{x})=\m{ax}$ for each $x\in J$.      
The toeplitz algebra associated to $(A,\al,L)$ is the universal $C^*$-algebra $\tp{A}{}$ generated by $A\cup
M$ with the relations of $A$, of $M$, the bi-module products and $m^*n=\langle m,n \rangle$. A redundancy in
$\tp{A}{L}$ is a pair $(a,k)\in A\times \widehat{K_1}$, ($\widehat{K_1}=\ovl{\text{span}}\{mn^*,m,n \in
M\}$), such that $am=km$ for every $m\in M$. The {\it Crossed Product by a Partial Endomorphism}
$\prcr{A}{}$ is the quotient of $\tp{A}{}$ by the ideal generated by all the elements $a-k$ where $(a,k)$ is 
a redundancy and $a\in \ker(\varphi)^{-1}\cap \varphi^{-1}(K(M))$.

In \cite{prodcruz} it was defined the algebra $\prcr{X}{}$. This algebra is constructed from a
partial $C^*$-dynamical
system $(C(X), \al, L)$ induced by a
local homeomorphism $\sigma:U\rightarrow X$, where $U$ is an open subset of a compact topological Hausdorff
space $X$. More specifically,
$$\funcao{\al}{C(X)}{C^b(U)}{f}{f\circ\sigma}$$ where $C^b(U)$ is the space of all continuous bounded functions
in $U$ and
$L:C_c(U)\rightarrow C(X)$ ($C_c(U)$ is the set of the continuous functions with compact support in $U$) is
defined by $$L(f)(x)=\left \{\begin{array}{cc} \sum\limits_{y\in
\sigma^{-1}(x)}f(y) & \text{if } x \in \sigma(U) \\
0 & \text{ otherwise}\end{array} \right .$$ for every $f\in C_c(U)$ and $x\in X$. 

In \cite{kwasniewski} it was defined the algebra $C^*(X,\al)$, called
covariance algebra. This algebra is also constructed from a partial dynamical system, that is, a continuous
map
$\sigma:U\rightarrow X$ where $X$ is a topologial compact Hausdorff space
, $U$ is a clopen subset of $X$ and $\sigma(U)$ is open.

If we suppose that $\sigma:U\rightarrow X$ is an local homeomorphism, $U$ clopen (and so $\sigma(U)$ is always
open) then $(X,\sigma)$ gives rise to two
 $C^*$-algebras, the covariance algebra $C^*(X,\sigma)$ and the crossed product by a partial endomorphism $\prcr{X}{}$.

In this paper we identify the transfer operators $L_\rho$ which may occur for $\al$. Moreover we
show that the constructions of covariance algebra and crossed product by partial
endomoprhism are not equivalent, in the following sense: we obtain examples of partial
dynamical systems $(X_A,\sigma_A)$ such that there does not exists invertible function
$\rho$ such that $\prcr{X_A}{\rho}=C^*(X,\al)$.

Given $C^*$-algebras $A$ and $B$, if we write $A=B$, we will say that $A$ and $B$ are
*-isomorphic.

\centerline{\large Acknowledgements}

The author wishes to express his thanks to R. Exel for the many stimulating conversations.

\section{Transfer operators of $X$ for $\al$}

Let $\sigma:U \rightarrow X$ be a local homeomorphism and $U$ an open subset of the compact Hausdorff space $X$.
This local homeomorphism induces the *-homomorphism  $$\funcao{\al}{C(X)}{C^b(U)}{f}{f\circ\sigma}.$$
Given a positive function $\rho\in C(U)$, for all $f\in C_c(U)$ we may define
$$L_\rho(f)(x)=\left \{\begin{array}{cc} \sum\limits_{y\in
\sigma^{-1}(x)}\rho(y)f(y) & \text{if } x \in \sigma(U) \\
0 & \text{ otherwise}\end{array} \right .$$ for each $x\in X$. Note that
$L_\rho(f)=L(\rho f)$, and since $\rho f\in C_c(U)$ and $L(\rho f)\in C(X)$ (see \cite{prodcruz}) then $L_\rho(f)$ in
fact is an element of $C(X)$.
In this way we may define the map
$L_\rho:C_c(U)\rightarrow C(X)$, which is linear and positive (by the fact that $\rho$
is positive). It is easy to see that $L_\rho(f\al(g))=L_\rho(f)g$ for each
$f\in C_c(U)$ and $g\in C(X)$.
The following proposition shows that if $U$ is clopen in $X$ then every transfer operator for $\al$ is of the
form $L_\rho$ for some $\rho\in C(U)$.

\begin{proposicao}
Let $L:C_c(U)\rightarrow C(X)$ ($U$ clopen in $X$) a transfer operator for $\al$, that is, $L$ is linear,
positive and
$L(g\al(f))=L(g)f$ for each $f\in C(X)$ and $g\in C_c(U)$. Then there exists $\rho\in C(U)$ such that $L=L_\rho$.
\end{proposicao}

\demo
Let $\{V_i\}_{i=1}^n$ be an open cover of $U$ such that $\sigma_{|_{V_i}}$ is a homeomorphism. (such cover
exists because $U$ is compact and $\sigma$ is a local homeomorphism). For each $i$ take an open subset $U_i\subseteq V_i$
such that $\ovl{U_i}\subseteq
V_i$ and $\{U_i\}_i$ is also a cover for $U$. Consider the partition of unity $\{\varphi_i\}_i$ subordinated to
$\{U_i\}_i$ and define $\xi_i=\sqrt{\varphi_i}$.
Since $\xi_i$ is positive for each $i$ then $L(\xi_i)$ is a positive function.
Define $\rho=\sum\limits_{i=1}^n\al (L(\xi_i))\xi_i$ which is also positive.
Given $f\in C_c(U)$ define for each $i$,

$$g_i(x)=\left\{\begin{array}{cc}
\xi_i(\sigma^{-1}(x))f(\sigma^{-1}(x)) & x \in \sigma(V_i) \\
0 & \text{otherwise}
\end{array} \right .$$

\af{ 1}{$g_i\in C(X)$ for all $i$}

Let $x_j\rightarrow x$. Suppose $x\in \sigma(V_i)$. Since $\sigma(V_i)$ is open we may suppose that $x_j\in
\sigma(V_i)$ for each $j$. Since $\sigma_{|_{V_i}}$ is a homeomorphism then $\sigma^{-1}(x_j)\rightarrow
\sigma^{-1}(x)$ in $V_i$ and so $g_i(x_j)=(\xi_i f)(\sigma^{-1}(x_j))\rightarrow (\xi_i
f)(\sigma^{-1}(x))=g_i(x)$. If $x\notin$ $\sigma(V_i)$ then $x\notin \sigma(\ovl{U_i})$, which is closed.
Therefore we may suppose that $x_j\notin \sigma(\ovl{U_i})$ and so $g_i(x_j)=0=g_i(x)$.
\fimaf

\af{ 2}{$\xi_i\al(g_i)=\varphi_i f$}

If $x\notin U_i$ then $(\xi_i\al(g_i))(x)=0=(\varphi_i f)(x)$. If $x\in U_i$ then
$\al(g_i)(x)=g_i(\sigma(x))=\xi_i(x)f(x)$ and so $\xi_i(x)\al(g_i)(x)=\xi^2(x)f(x)=\varphi(x)f(x)$.
\fimaf

Since $\varphi$ is partition of unity then $f=\sum\limits_{i=1}^n\varphi_i
f=\sum\limits_{i=1}^n\xi_i\al(g_i)$, where the last equality follows by claim 2. Then

$$L(f)=\sum\limits_{i=1}^nL(\xi_i\al(g_i))=\sum\limits_{i=1}^nL(\xi_i)g_i.$$

We show that $L=L_\rho$. If $x\notin \sigma(U)$ then $L_\rho(f)(x)=0=L(f)(x)$ by definition.

Given $x\in \sigma(U)$,

$$L_\rho(f)(x)=\sum\limits_{y\in \sigma^{-1}(x)}\rho(y)f(y)=\sum\limits_{y\in
\sigma^{-1}(x)}\sum\limits_{i=1}^n\al(L(\xi_i))(y)\xi_i(y)f(y)=$$

$$=\sum\limits_{y\in\sigma^{-1}(x)}\sum\limits_{i:y\in U_i}L(\xi_i)(x)\xi_i(y)f(y).$$

On the other hand,

$$L(f)(x)=\sum\limits_{i=1}^nL(\xi_i)(x)g_i(x)=\sum\limits_{i:x\in \sigma(U_i)}
L(\xi_i)(x)\xi_i(\sigma^{-1}(x))f(\sigma^{-1}(x)).$$ To see that those two expressions are equal note that the
summands are the same. \fim

Denote by $M_\rho$ the Hilbert bi-module generated by $C_c(U)$ with the inner product given by $L_\rho$ and by
$\widehat{K_{1\rho}}$ the algebra generated by $nm^*$ in $\tp{X}{\rho}$. Moreover, denote by
$\varphi_\rho:C(X)\rightarrow L(M_\rho)$ the *-homomorphism given by the left product of $A$ by $M_\rho$.

\begin{lema}
Let $\rho, \rho'\in C(U)$ positive functions. If $\ker(\rho)=\ker(\rho')$ then
$\ker(\varphi_\rho)=\ker(\varphi_{\rho'})$.
\end{lema}

\demo Let $f\in C(X)$.
Then $f\in \ker(\varphi_\rho)\Leftrightarrow fm=0$ for each $m\in M_\rho$ $\Leftrightarrow$ $\m{fg}=f\m{g}=0$
for each $g\in C_c(U)$.
It is easy to check that $\m{fg}=0$ in $M_\rho$ if and only if $\rho fg=0$.
Then $f\in \ker(\varphi_\rho)$ if and only if $\rho fg=0$ for each $g\in C_c(U)$. 
In the same way, $f\in \ker(\varphi_\rho')$ if and only if $\rho' fg=0$ for each $g\in C_c(U)$. 
Since $\ker(\rho)=\ker(\rho')$ then $\rho fg=0$ if and only if $\rho' fg=0$ for each $g\in C_c(U)$.
\fim

\begin{proposicao}
If $\rho$ and $\rho'$ are elements of $C(U)$ such that there exists $r\in C(U)$ such that $r(x)\neq 0$ for
each $x\in
U$ and $\rho=r\rho'$ then $\prcr{X}{\rho}$ and $\prcr{X}{\rho'}$
are *-isomorphic.
\end{proposicao}

\demo
Let us define a *-homomorphism from $\prcr{X}{\rho}$ to $\prcr{X}{\rho'}$. Define
$$\funcao{\psi_1}{C(X)}{\tp{X}{\rho'}}{f}{f}.$$ Let $\xi=\sqrt{r}$, and note that for each
$g\in C_c(U)$,
$$\|\m{g}\|_\rho^2=\|L_\rho(g^*g)\|=\|L(\rho g^*g)\|=\|L(r\rho' g^* g)\|=\|L_{\rho'}((\xi g)^*\xi
g)\|=\|\m{\xi g}\|_{\rho'}^2.$$ So we may define $\psi_2:M_\rho\rightarrow
\tp{X}{\rho'}$ by
$\psi_2(\m{g})=\m{\xi g}$. Let $\psi_3=\psi_1\cup \psi_2$. We show that $\psi_3$ extends
to $\tp{X}{\rho}$. For each
$f\in C(X)$ and $g\in
C_c(U)$ we have $$\psi_3(f)\psi_3(\m{g})=f\m{\xi g}=\m{\xi fg}=\psi_3(\m{fg})$$ and
$$\psi_3(\m{g})\psi_3(f)=\m{\xi g}f=\m{\xi g \al(f)}=\psi_3(\m{g\al f}).$$ Moreover, if $h\in C_c(U)$ then
$$\psi_3(\m{g})^*\psi_3(\m{h})=\m{\xi g}^*\m{\xi h}=L_{\rho'}((\xi g)^*\xi
h)=L_{\rho'r}(g^*h)=L_\rho(g^*h)=\psi_3(L_\rho(g^*h)).$$ So $\psi_3$ axtends to $\tp{X}{\rho}$. Let $(f, k)\in
C(X)\times \widehat{K_{1\rho}}$ a
redundancy with \lb $f\in \ker(\varphi_\rho)^\bot \cap \varphi_\rho^{-1}(K(M_\rho))$. Since
$\psi_3(M_\rho)\subseteq M_{\rho'}$ it follows that $\psi_3(k)\in \widehat{K_{1\rho'}}$ and so $(\psi_3(f), \psi_3(k))\in
C(X)\times \widehat{K_{1\rho'}}$. Moreover, given $g\in C_c(U)$ then $\xi^{-1}g\in C_c(U)$ and
$\psi_3(\m{\xi^{-1}g})=\m{g}$ from where $\psi_3(M_\rho)$ is dense in $M_{\rho'}$, and so, since $fm=km$ for each $m\in
M_\rho$ then $\psi_3(f)n=\psi_3(k)n$ for every $n\in M_{\rho'}$. Therefore $(\psi_3(f), \psi_3(k))$ is a
redundancy. Since $f\in \ker(\varphi_\rho)^\bot$, by the previous lemma, $\psi_3(f)\in\ker(\varphi_{\rho'})^\bot$.
Then, since
$(\psi_3(f),\psi_3(k))$ is a redundancy of $\tp{X}{}$ then by  [\ref{prodcruz}:2.6], $\psi_3(f)\in
\varphi^{-1}(K(M_{\rho'}))$.
So $\psi_3(f)\in \ker(\varphi_{\rho'})^\bot \cap \varphi_{\rho'}^{-1}(K(M_{\rho'}))$. This shows that if $\phi$ is the quotient
*-homomorphism from
$\tp{X}{}$ in $\prcr{X}{}$ then $\phi\circ\psi_3:\tp{X}{\rho}\rightarrow \prcr{X}{}$ is a homomorphism which
vanishes on all the elements of the form $(a-k)$ where $(a,k)$ is a redundancy and
$a\in\varphi_\rho^{-1}(K(M_\rho))\cap
\ker(\varphi_\rho)^\bot$. So we obtain a *-homomorphism
$$\dfuncao{\psi}{\prcr{X}{\rho}}{\prcr{X}{\rho'}}{f}{f}{\m{g}}{\m{\xi g}}.$$

In the same way we may define the *-homomorphism
$$\dfuncao{\psi_0}{\prcr{X}{\rho'}}{\prcr{X}{\rho}}{f}{f}{\m{g}}{\m{\xi^{-1}g}}.$$ Note that $\psi_0$ is the
inverse of
$\psi$, showing that tha algebras are *-isomorphic.\fim

\begin{corolario}\label{c1}
If $\rho\in C(U)$ is a positive function such that $\rho(x)\neq 0$ for all $x\in U$ then $\prcr{X}{\rho}$ is
*-isomorphic to $\prcr{X}{}$.
\end{corolario}

\demo Note that the transfer operator $L$ associated to the algebra $\prcr{X}{}$ is the
operator $L_{1_U}$. Since $\rho=1_U$ is invertible, taking $r=\rho^{-1}$, by the previous
proposition follows the corollary.\fim

\section{Relationship between the Covariance Algebra and the Crossed Product by Partial Endomorphism}

We show here that given a partial dynamical system $\sigma:U\rightarrow X$, where $U$ is
clopen, there exists an other partial dynamical system $\m{\sigma}:\m{U} \rightarrow\m{X}$ (called
in \cite{kwasniewski} the $\sigma$-extension of $X$) such that
$C^*(\sigma,X)=\prcr{\m{X}}{}$. Moreover, if $\sigma$ is injective then $C^*(\sigma,X)=\prcr{X}{}$.

\subsection{The Covariance Algebra as an Crossed Product by a Partial Endomorphism}

Let us start with a summary of the construction of the covariance algebra. Let
$\sigma:U\rightarrow X$ a continuous map, $U\subseteq X$ clopen, $X$ compact Hausdorff and
$\sigma(U)$ open. Denote $\sigma(U)=U_{-1}$. Consider the space $X\cup \{0\}$, where
$\{0\}$ is a symbol, which we define to be clopen. So $X\cup\{0\}$ is a compact and
Hausdorff space.

Define $\m{X}\subset \prod\limits_{i=0}^\infty X\cup \{0\}$, $$\m{X}=\bigcup\limits_{N=0}^\infty X_N\cup X_\infty$$ onde
$$X_N=\{(x_0,x_1,...,x_N,0,0,...): \sigma(x_i)=x_{i-1} \text{ e } x_N\notin U_{-1}\}$$ and
$$X_\infty=\{(x_0,x_1,...):\sigma(x_i)=x_{i-1}\}.$$ In $\m{X}$ we consider the product topology induced from
$\prod X\cup \{0\}$.

By $[\ref{kwasniewski}: 2.2]$ $\m{X}$ is compact. Define
$$\funcao{\Phi}{\m{X}}{X}{(x_0,x_1,...)}{x_0}$$ which is continuous and surjective. Consider
the clopen subsets $\m{U}=\Phi^{-1}(U)$ and
$\m{U_{-1}}=\Phi^{-1}(U_{-1})$ and the continuous map
$$\funcao{\m{\sigma}}{\m{U}}{\m{U_{-1}}}{(x_0,x_1,...)}{(\sigma(x_0),x_0,x_1,...)}.$$ Those maps satisfies the
relation $$\Phi(\m{\sigma}(\m{x}))=\sigma(\Phi(\m{x})).$$

Note that $\m{\sigma}$ is in fact an homeomorphism. This homeomorphism induces the
*-isomorphism
$$\funcao{\theta}{C(\m{U_{-1}})}{C(\m{U})}{f}{f\circ \m{\sigma}}.$$ So we may consider the partial crossed
product
$C(\m{X})\rtimes_\theta \Z$ (see \cite{exelprodcruzpartial}).

\begin{definicao}[\ref{kwasniewski}: 4.2]
The covariance algebra associated to the partial system $(X,\sigma)$ is the algebra $C(\m{X})\rtimes_\theta
\Z$ and will be denoted $C^*(X,\sigma)$.
\end{definicao}

\begin{lema}\label{l1}
If $\sigma:U\rightarrow X$ is injective, $U$ clopen and $U_{-1}$ open then
$C(X)\rtimes_{\theta}\Z=\prcr{X}{}$, where $\theta:C(U_{-1})\rightarrow C(U)$
is given by $\theta(f)=f\circ\sigma$.
\end{lema}

\demo
Define $\psi_1:C(X)\cup M\rightarrow C(X)\rtimes_{\theta}\Z$ by $\psi_1(f)=f\delta_{0}$ and
$\psi_1(\m{1_U})=1_U\delta_1$. It is easy to check that $\psi_1$ axtends to $\tp{X}{}$. We show that
$\Psi_1$ vanishes on the redundancies. Let $(f,k)$ redundancy with $f\in \ker(\varphi)^\bot\cap
\varphi^{-1}(K(M))$. By [\ref{prodcruz}: 2.6], $f\in C(U)$. Then
$\psi_1(\m{f})\psi_1(\m{1_U})^*=f\delta_11_{U_{-1}}\delta_{-1}=\theta(\theta^{-1}(f)1_{U_{-1}})\delta_0=\psi_1(f)$.
Take $(k_n)_n\subseteq \widehat{K_1}$, $k_n=\sum\limits_{i}m_{ni}l_{ni}^*$ where $m_{ni},l_{ni}\in M$. Then
$$(\psi_1(f)-\psi(k))(\psi_1(f)-\psi(k))^*=(\psi_1(f)-\psi_1(k))\psi_1(f-k)=\psi_1(f-k)(\psi_1(\m{f}\m{1_U}^*)-\psi_1(k))^*=$$
$$=\psi(f-k)(\m{1_U}\m{f}^*-k)=\lim\limits_{n\rightarrow \infty}\psi(f-k)(\m{1_U}\m{f}^*-k_n)=0.$$
The last equality follows by the fact that $(f-k)m=0$ for each $m\in M$.
So, by passage to the quotient we may consider $\psi:\prcr{X}{}\rightarrow C(X)\rtimes_{\theta}\Z$. By the
other hand, define
$$\funcao{\psi_0}{C(X)}{\prcr{X}{}}{f}{f}$$ which is a *-homomorphism. Note that for each $f\in C(U_{-1})$,
$$\m{1_U}\psi_0(f)\m{1_U}^*=\m{1_U\al(f)}\m{1_U}^*=1_U\al(f)=\theta(f)=\psi_0(\theta(f))$$ and moreower
$\m{1_U}$ is a partial isometry such that $\m{1_U}\m{1_U}^*=1_{U}$ and $\m{1_U}^*\m{1_U}=1_{U_{-1}}$.
Then, since $(\psi_0,\m{1_U})$ is a covariant representation of $C(X)$ in $\prcr{X}{}$, there exists a
*-homomorphism $\psi':C(X)\rtimes_{\theta}\Z\rightarrow \prcr{X}{}$ such that
$\psi'(f\delta_n)=f\m{1_U}^n$ (see [\ref{exelprodcruzpartial}: 5]).
The *-homomorphisms $\psi$ and $\psi'$ are inverses of each other, and so the algebras are *-isomorphic.\fim

\begin{corolario}
$C^*(X,\sigma)=\prcr{\m{X}}{}$
\end{corolario}

\demo
Follows by the definition of covariance algebras and by the previous lemma.
\fim

By the following proposition, if $\sigma$ is injective then the constructions of
covariance algebra and crossed product by partial endomorphism are equivalent.

\begin{proposicao}
If $\sigma:U\rightarrow X$ is injective then $C^*(X,\sigma)=\prcr{X}{}$.
\end{proposicao}

\demo By [\ref{kwasniewski}: 2.3] the map $$\funcao{\Phi}{\m{X}}{X}{(x_0,x_1,...)}{x_0}$$ is a homeomorphism.
Moreover, since $\Phi\circ \m{\sigma}=\sigma\circ\Phi$ then
$C(\m{X})\rtimes_{\m{\theta}}\Z=C(X)\rtimes_{\theta}\Z$. By the previous lemma $C(X)\rtimes_{\theta}\Z=\prcr{X}{}$.
\fim

\subsection{Cuntz-Krieger algebras}

We show examples of partial dynamical system $\sigma_A:U\rightarrow X_A$ such that there
does not exists an invertible function $\rho\in C(U)$ such that $\prcr{X}{\rho}$ and
$C^*(X,\al)$ are *-isomorphic. The examples are based on the Cuntz-Krieger algebras.

Let $A$ be a $n\times n$ matrix with $a_{i,j}\in \{0,1\}$. Denote by $Gr(A)$ the directed
graph of $A$, that is, the vertex set is $\{1,...,n\}$ and $A(i,j)$ is the number of oriented
edges from $i$ to $j$. A path is a sequence $x_1,...,x_m$ such that $A(x_i,x_{i+1})=1$ for
each $i$. The graph $Gr(A)$ is transitive if for each $i$ and $j$ there exists a path
from $i$ to $j$, that is, a path $x_1,...,x_m$ such that $x_1=i$ and $x_m=j$. The graph
is a cicle if for each $i$ there exists only one $j$ such that $A(i,j)=1$.

Let
$$X_A=\{x=(x_1,x_2,..)\in \{1,...n\}^\N\,\,:\,\,A(x_i,x_{i+1})=1 \,\,\forall \,\, i\}\subseteq \{1,...,n\}^\N$$ and
$$\funcao{\sigma_A}{X_A}{X_A}{(x_0,x_1,...)}{(x_1,x_2,...)}.$$ Consider the set
$$\ovl{X_A}=\{(x_i)_{i\in \Z}\in \{1,...,n\}^\Z \,\,:\,\, A(x_i,x_i+1)=1\,\, \forall \,\,i\}\subseteq
\{1,...,n\}^\Z$$
 and the map $\ovl{\sigma_A}:\ovl{X_A}\rightarrow \ovl{X_A}$ defined by $\ovl{\sigma_A}((x_i)_{i\in
 Z})=(x_{i+1})_{i\in \Z}$. It is showed in [\ref{kwasniewski}: 2.8] that there exists a homeomorphism $\Phi:\m{X_A}\rightarrow \ovl{X_A}$
 such that $\Phi\circ\m{\sigma_A}=\ovl{\sigma_A}\circ \Phi$. Therefore
 $\prcr{\m{X_A}}{}=\prcr{\ovl{X_A}}{}$ and so $C^*(X_A,\sigma_A)=\prcr{\ovl{X_A}}{}$.
So we may analyze the ideal structure of $C^*(X_A,\sigma_A)$ by using the theory
developed for $\prcr{\ovl{X_A}}{}$ in \cite{prodcruz}. This theory is based on the
$\ovl{\sigma_A}$,$\ovl{\sigma_A}^{-1}$ invariant open subsets of $\ovl{X_A}$. (In a system $\sigma:U\rightarrow
X$, a subset $V\subseteq X$ is $\sigma,\sigma^{-1}$ invariant if $\sigma(U\cap V)\subseteq V$ and
$\sigma^{-1}(V)\subseteq V$).

\begin{proposicao}\label{abertos invariantes}
If $Gr(A)$ is transitive and is not a cicle then there exists at least one open non
trivial $\ovl{\sigma_A}$, $\ovl{\sigma_A}^{-1}$ invariant subset of $\ovl{X_A}$.
\end{proposicao}

\demo Let $r=x_1,x_2,...,x_n$ an admissible word (that is, $A(x_i,x_{i+1})=1$ for each $i$). Let
$V_r=\{x\in \ovl{X_A}:r\in x\}$. Note that $V_r$ is open and $\ovl{\sigma_A}$, $\ovl{\sigma_A}^{-1}$ invariant.
We show that there exists a such non trivial $V_r$.
Take $x_1\in \{1,...,n\}$. Consider an admissible word $x_1,...,x_m$ where $x_j\neq x_1$ for each $j>1$ and
$A(x_m,x_1)=1$. Such word exists because $Gr(A)$ is transitive. Let $r=x_1,...,x_m,x_1$. Then
$$y=(...,x_m,\stackrel{\circ}{x_1},x_2,...,x_m,x_1,x_2,...)\in V_r$$ where $\stackrel{\circ}{x_1}$ means 
$y_0=x_1$. 

We conclude the proof by showing that $V_r\neq \ovl{X_A}$.
Suppose that there exists $y_0\in \{1,...,n\}$ with $y_0\notin\{x_1,...,x_m\}$. Let
$x_1,y_1,...,y_t,y_0,s_1,...,s_l$
an admissible word
such that $y_j\neq x_1$ and $s_j\neq x_1$ for each $j$ and $A(s_l,x_1)=1$. Then
$$(...s_l,\stackrel{\circ}{x_1},y_1,...,y_t,y_0,s_1,...,s_l,x_1...)\notin V_r.$$ If
$\{x_1,...,x_m\}=\{1,...,n\}$, since $Gr(A)$ is not a cicle, for some $x_i$ there exists $x_t$ such that
$A(x_i,x_t)=1$ and $x_t\neq x_{i+1}$. (if $i=m$ consider $x_{i+1}=x_1$). If $x_t=x_1$ (and so $i\neq m$)
consider an admissible word $x_1,...,x_i,x_1$ and note that
$$(...,x_i,\stackrel{\circ}{x_1},x_2,...,x_i,x_1,...)\notin V_r.$$ If $x_t\neq x_1$ consider an admissible word
$x_1,x_2,...,x_i,x_t,y_1,...,y_l$ such that $y_j\neq x_1$ and $A(y_l,x_1)=1$ (if there does not exists $y_1\neq x_1$
such that $A(x_t,y_1)=1$ then $y_1,...,y_l$ is the empty word) and so
$$(...,y_l,\stackrel{\circ}{x_1},x_2,...x_i,x_t,y_1,...,y_l,y_1,...)\notin V_r.$$
So $V_r\neq \ovl{X_A}$.\fim

Now we analyse the $\sigma_A,\sigma_A^{-1}$ invariant subsets of $X_A$.

\begin{proposicao}\label{p2}
If $Gr(A)$ is transitive and is not a cicle the the unique open $\sigma_A$-invariant subset of $X_A$ are
$\emptyset$ and $X_A$.
\end{proposicao}

\demo Let $V\subseteq X_A$ an open nonempty $\sigma_A$ invariant susbet of $X_A$. Let $x\in V$ and $V_m$
a open neighbourhood of $x$, $V_m\subseteq V$, $$V_m=\{y\in X_A: x_i=y_i \text{ for each } 1\leq i\leq m\}.$$
Given $z\in X_A$ take
$r=r_1,...,r_t$ a path from $x_m$ to $z_1$. Then $$s=(x_1,...,x_m, r_2,...,r_{t-1},z_1,z_2,...)\in V_m$$ and since
$V$ is $\sigma_A$ invariant then $z=\sigma_A^{m+t-2}(s)\in V$.
So $V=X_A$. \fim

According \cite{prodcruz} a partial dynamical system $\sigma:U\rightarrow X$ is topologically free if the
closure of $V^{i,j}=\{x\in U:\sigma^i(x)=\sigma^j(x)\}$ has empty interior for each $i,j\in \N$, $i\neq j$.

\begin{proposicao}\label{p3}
If $Gr(A)$ is transitive and is not a cicle then $(X_A, \sigma_A)$ is topologically free.
\end{proposicao}

\demo Suppose that $\ovl{V^{i,j}}$ has nonempty interior and $i<j$, $j=i+k$. Let $x'$ be an interior point of
$\ovl{V^{i,j}}$ and $V_{x'}\subseteq\ovl{V^{i,j}}$ open neighbourhood of $x'$. Take $x\in V^{i,j}\cap V_{x'}$
. Since
$\sigma_A^i(x)=\sigma_A^j(x)$ then $z_{i+t}=z_{j+t}$ for each $t\in \N$ and since $j=i+k$ then
$x=(x_1,...,x_{i-1},r,r,...)$ where $r=x_ix_{i+1}...x_{i+k-1}$. Consider the open subset $$V_m=\{z\in
X_A:z_i=x_i, 1\leq i\leq m\}$$ where $m$ is such that $m\geq i+k$ and $V_m\subseteq V_{x'}$ . Then, if
$y\in V_m$ with $y\in V^{i,j}$ then $y=x$. Therefore $V_m=\{x\}$. 
We show that there exists $z\in V_m$ with $z\neq
x$, and that will be a contradiction. Suppose $y_0\in \{1,...,n\}$ and $y_0\notin \{x_i,...,x_{i+k-1}\}$. Take 
a path $s=s_1,...,s_t$
from $x_i$ to $x_{i+k-1}$ such that $s_j=y_0$ for some $j$. Then $z=(x_1,...,x_{i-1},r,r,...,r,s,s,...)\in
V_m$ (where $r$ is repeated $m$ times) but $z\neq x$. Suppose $\{1,...,n\}=\{x_i,...,x_{i+k-1}\}$. Since
$Gr(A)$ is not a cicle then for some $x_j$ there exists $x_t\neq
x_{j+1}$ (consider $x_{j+1}=x_i$ if $j=i+k-1$) such that $A(x_j,x_t)=1$. Let $s$ be a path from $x_t$ to
$x_{i+k-1}$ and define $p=x_i,...,x_j,x_t,s$.
Then $$z=(x_1,...,x_{i-1},r,r,...,r,x_i,...,x_j,x_t,p,p,p,,,...)\in V_m$$ (where $p$ is repeated $m$ times) but $z\neq x$.
So, it is showed that there exists $z\in V_m$, $z\neq x$. Therefore, $\ovl{V^{i,j}}$ has empty interior for
each $i,j$.\fim

\begin{teorema}
If $Gr(A)$ is transitive  and is not a cicle then $C^*(X_A,\sigma_A)$ and
$\prcr{X_A}{}$ are not *-isomorphic $C^*$-algebras.
\end{teorema}

\demo By $\ref{l1}$, $C^*(X_A,\sigma_A)=\prcr{\m{X_A}}{}$ and since $\prcr{\m{X_A}}{}=\prcr{\ovl{X_A}}{}$ then
\lb $C^*(X_A,\sigma_A)=\prcr{\ovl{X_A}}{}$. By $\ref{abertos invariantes}$, $\ovl{X_A}$ has at least one non
trivial open $\ovl{\sigma_A},\ovl{\sigma_A}^{-1}$ invariant subset and by
[\ref{prodcruz}: 3.9] $\prcr{\ovl{X_A}}{}$ has
at least on non trivial ideal. On the other hand, by $\ref{p2}$, $(X_A,\sigma_A)$ has no open
$\sigma_A,\sigma_A^{-1}$ invariant subsets and by $\ref{p3}$, $(X_A, \sigma_A)$ is topologically free. By
[\ref{prodcruz}: 4.8] $\prcr{X_A}{}$ is simple.
So $C^*(X_A,\sigma_A)$ and $\prcr{X_A}{}$ are not *-isomorphic. \fim

\begin{corolario}
If $Gr(A)$ is transitive and is not a cicle then there does not exists transfer operator $L_\rho$, with $\rho(x)\neq
0$ for each $x\in U$ such that $C^*(X_A,\sigma_A)$ and $\prcr{X_A}{}$ are *-isomorphic $C^*$-algebras.
\end{corolario}

\demo Follows by the previous theorem an by \ref{c1}.\fim

\vspace{2cm}

Departamento de Matemática, Universidade Federal de Santa Catarina, Brasil.\nl
E-mail: royer@mtm.ufsc.br

\end{document}